\documentclass{amsart}

\usepackage[latin1]{inputenc}
\usepackage{amssymb}
\usepackage{latexsym}
\usepackage{longtable}

\newtheorem{theorem}{Theorem}[section]

\newtheorem{lemma}[theorem]{Lemma}
\newtheorem{proposition}[theorem]{Proposition}

\newtheorem{remark}{Remark}[section]

\newcommand\Z{{\mathbb Z}}
\newcommand\R{{\mathbb R}}
\newcommand\F{{\mathbb F}}
\newcommand\Q{{\mathbb Q}}
\newcommand\OO{{\mathcal O}}
\newcommand\cR{{\mathcal R}}
\newcommand\T{{\mathcal T}}
\newcommand\M{{\mathrm {M}}}
\newcommand\pp{{\mathfrak p}}
\newcommand\nbar{\nmid}
\renewcommand{\mod}[1]{\; (\mbox{\rm mod } #1)}
\newcommand\Gal{\mbox{\rm Gal}}
\newcommand\End{\mbox{\rm End}}

\newcommand\Tr{\mbox{\rm Tr}}
\newcommand\n{\mbox{\rm n}}
\newcommand\Pic{\mbox{\rm Pic}}

\begin{document}

\pagestyle{plain}

\title{Easy decision-Diffie-Hellman groups}

\author{Steven D. Galbraith
and Victor Rotger}
\footnote{The first author thanks the Nuffield foundation
grant NUF-NAL-02 for support.}

\address{ Mathematics Department, Royal Holloway University of London,
Egham, Surrey TW20 0EX, UK. \\
Universitat Polit\`{e}cnica de Catalunya,
Departament de Matem\`{a}tica Aplicada IV (EUPVG), Av.\ Victor
Balaguer s/n, 08800 Vilanova i la Geltr\'{u}, Spain.}
\email{Steven.Galbraith@rhul.ac.uk, vrotger@mat.upc.es}

\begin{abstract}
The decision-Diffie-Hellman problem (DDH) is
a central computational problem
in cryptography.
It is known that the Weil and Tate pairings
can be used to solve many DDH
problems on elliptic curves.
Distortion maps are an
important tool for solving DDH problems using pairings
and it is known that distortion maps exist for all
supersingular elliptic curves.
We present an algorithm to
construct suitable distortion maps.
The algorithm is efficient on the curves usable in
practice, and hence all DDH problems on
these curves are easy.
%
We also discuss the issue of which DDH problems on ordinary
curves are easy.
\end{abstract}

\maketitle

\section{Introduction}

It is well-known that the Weil and Tate pairings make many
decision-Diffie-Hellman (DDH) problems on elliptic curves
easy.  This observation is behind exciting new
developments in pairing-based cryptography.
This paper studies the question of which
DDH problems are easy and which are not necessarily
easy.
First we recall some definitions.

\vskip 0.2cm

\noindent {\bf Decision Diffie-Hellman problem (DDH):}
Let $G$ be a cyclic group of
prime order $r$ written additively.  The DDH problem is to
distinguish the two distributions in $G^4$
\begin{eqnarray*}
    D_1 &=& \{ ( P, aP, bP, abP) : P \in G, 0 \le a, b < r \} \;\;
  \mbox{ and}\\
    D_2 &=& \{ ( P, aP, bP, cP) : P \in G, 0 \le a, b, c < r \}.
\end{eqnarray*}
Here $D_1$ is the set of valid Diffie-Hellman-tuples and $D_2 = G^4$.
By `distinguish' we mean there is an algorithm which takes
as input an element of $G^4$ and outputs ``valid'' or
``invalid'', such that if the input is chosen with
probability 1/2 from each of $D_1$ and $D_2 - D_1$ then the
output is correct with probability
significantly more than 1/2.
(for precise definitions see Boneh~\cite{boneh}).
The DDH problem for a family of groups is
said to be hard if there is  no polynomial time algorithm
which distinguishes the two distributions.
A widely believed assumption in cryptography is that
there exist families of groups for which the DDH is problem
is hard.

We now give a generalisation of the DDH problem
which, following Boneh, Lynn and Shacham~\cite{bls},
we call co-DDH.

\vskip 0.2cm

\noindent {\bf Generalised Decision Diffie-Hellman problem
(co-DDH):}
Let $G_1$ and $G_2$ be two cyclic groups of prime order $r$.
The co-DDH problem is to
distinguish the two distributions in $G_1^2 \times G_2^2$
\[
  \begin{array}{l}
    \{ ( P, aP, Q, aQ) : P \in G_1, Q \in G_2, 0 \le a < r \} \; \;
  \mbox{ and } \\
    \{ ( P, aP, Q, cQ) : P \in G_1, Q \in G_2, 0 \le a, c < r \}.
  \end{array}
\]

The goal of this article is to determine which
DDH and co-DDH problems on elliptic curves are made easy by
using pairings.
A common technique is to use distortion maps
(endomorphisms which map certain subgroups of $E[r]$
to different subgroups)
to ensure that the required pairing values are non-trivial.
Theorem 5 of Verheul~\cite{ver-full} states that
a suitable distortion map always exists for
subgroups of supersingular curves.
This result alone does not imply that all
DDH problems can be solved efficiently, since we
require an explicit description of the map.

In Sections~\ref{sec-2} and \ref{sec-3}
we show that the trace map handles almost all cases.
In Section~\ref{main-theorem} we
give an alternative proof of Theorem 5 of \cite{ver-full}
(restricting to the remaining cases)
which is more constructive.
In Section~\ref{algorithm} we show that a certain
endomorphism $\sqrt{-d}$ suffices, and we
give an algorithm to construct this distortion map.
The complexity analysis of our algorithm proves that
all DDH problems are easy on the supersingular elliptic
curves which could potentially be used in practice.
Sections~\ref{sec-7} and \ref{not-iso} illustrate
the theory in concrete situations.
In particular, Section~\ref{sec-7} lists some well-known examples
and shows that they always suffice in practice.
Section~\ref{not-iso} gives examples of our
method in the case where the distortion map
cannot be an automorphism of the curve.
Some of the examples in Sections~\ref{sec-7} and \ref{not-iso}
show that our algorithm is not optimal, in the sense that
it does not necessarily produce an endomorphism of minimal degree.

Our results may have applications
as they mean that
cryptographic protocols can use random
points $P, Q$ on a supersingular elliptic curve and there is always
a modified pairing so that $e( P, Q ) \ne 1$.

In the case of ordinary elliptic curves
there are two hard DDH subgroups remaining.
Understanding whether these are truly hard
is a challenge to any interested person.

\section{Elliptic curves}
\label{sec-2}

We will be concerned with elliptic curves $E$ over
finite fields $\F_q$ such that $r$ is a large prime
dividing $\# E( \F_q )$ and such that $\gcd(r,q)=1$.
The embedding degree is the smallest positive integer
$k$ such that $r \mid (q^k - 1)$.
We restrict attention to elliptic curves such that
$k$ is not large (say, bounded by a fixed polynomial in $\log(q)$).
Hence, one can efficiently compute in $E( \F_{q^k} )$.
We always assume that $k$ is coprime to $r$
(this is always true since $r$ is a large prime and $k$ is small).

We will repeatedly make use of the following properties
of the Weil pairing (see Silverman \cite{silv} Section III.8).

\begin{lemma}
\label{lemma-1}
Let $E$ be an elliptic curve over $\F_q$ and let $P, Q \in E( \F_q)$
be points of prime order $r$. Then
\begin{enumerate}
\item $e_r( P, P ) = 1$.
\item If $P$ and $Q$ generate
$E[r]$ then $e_r( P, Q ) \ne 1$.
\item $R \in \langle P \rangle $ if and only if $e_r( P, R ) = 1$.
\end{enumerate}
\end{lemma}

\begin{proof}
The first statement is the well-known alternating property
of the Weil pairing.

Property 2 follows since if $e_r(P, Q)=1$
then $e_r( P, aP+bQ) = 1$ for all $a, b \in \Z$
which contradicts non-degeneracy
of the Weil pairing.

If $R = aP + bQ$ then $e_r( P, R) = e_r( P, Q)^b$
and this is 1 if and only if $b \equiv 0 \mod{r}$.
This proves property 3.
\end{proof}

\begin{remark}
Property 3 shows that
the subgroup membership problem for any cyclic
subgroup $G \subset E( \F_q )$
is easily solved using the Weil pairing if the
embedding degree is small.
Note that property 3 does not necessarily hold for
the Tate pairing (for details on the Tate pairing
see Frey and R\"uck~\cite{fr} or \cite{gal}).
\end{remark}

The above properties clearly imply that all genuine co-DDH problems
are easy.  This result is already well-known,
but for emphasis we state it as a proposition.

\begin{proposition}
\label{main-prop}
Let $E$ be an elliptic curve over $\F_q$ and let $r$
be a prime.  Suppose that $E[r] \subset E( \F_{q^k} )$
where $k$ is polynomial in $\log(q)$.
Let $G_1 , G_2 $ be cyclic subgroups of order $r$ in $E( \F_{q^k} )$
such that $G_1 \ne G_2$.  Then all co-DDH problems in
$G_1, G_2$ can be
solved in polynomial time.
\end{proposition}

\begin{proof}
The fact $G_1 \ne G_2$ implies $G_1 \cap G_2 = \{ 0_E \}$.
Hence for all $P \in G_1$, $Q \in G_2$, with $ P, Q \ne 0_E$
we have $\{ P, Q \}$ forming a basis for $E[r]$ and so by
property (2), $e_r( P, Q ) \ne 1$.

The co-DDH problem on a tuple $(P_1, P_2, Q_1, Q_2)$
is therefore solved by testing whether
\[
    e_r( P_1, Q_2) \stackrel{?}{=} e_r( P_2, Q_1 ).
\]
\end{proof}

\begin{remark}
As mentioned above, this result is not always
true for the Tate pairing.
However, in most practical cases the Tate pairing can be
used, and will give a
more efficient solution (see~\cite{bkls, ghs, gal}
for details).
\end{remark}

For the remainder of the paper we will be concerned
with solving DDH problems.  Clearly, the Weil pairing
cannot directly be used to solve these problems.

When $k=1$ and
$E( \F_q )[r]$ is a cyclic group of order $r$ then,
due to the non-degeneracy of the Tate pairing,
the DDH problem in this group can be solved in polynomial time.

The case $k=1$ and $E( \F_q )[r]$ non-cyclic is more interesting
(the curve $E$ is ordinary whenever $r$ is large).
This case has been considered by Joux and Nguyen~\cite{jn}.
The Weil and Tate pairings can have very different behaviour in
this case (for example, there are cases where the Tate
pairing always gives non-trivial self
pairings and cases where the Tate pairing never gives non-trivial
self pairings).
Theorem 7 of \cite{ver-full} shows that many DDH problems can be
solved in this case (using the Weil pairing
with a suitable distortion map).

In practice, the case $k>1$ is of greater interest.  Hence,
for the remainder of the paper we make the following assumption:

$\\ $ {\bf The embedding degree is assumed to be $k \ge 2$.}

\section{Trace maps}
\label{sec-3}

The trace map was proposed as a distortion map by Boneh et al
in the full versions of
\cite{bf} and \cite{bls}.
Since $\F_{q^k} / \F_q$ is a Galois extension
we can define, for any point $P \in E( \F_{q^k} )$,
\[
     \Tr( P ) = \sum_{i=0}^{k-1} \pi^i( P )
\]
where $\pi$ is the $q$-power Frobenius map.
Equivalently, if $P=(x,y)$ then
\[
     \Tr( P ) = \sum_{i=0}^{k-1} ( x^{q^i}, y^{q^i} ).
\]
The trace map is a group homomorphism and
if $P \in E( \F_q )$ then $\Tr(P) = kP$.

Let $P, Q \in E( \F_{q^k})[r]$.
Define the function $e( P, Q )$ to be either the
Weil pairing $e(P,Q) = e_r( P, Q )$
or the Tate pairing $e(P,Q) = \langle P, Q \rangle_r^{(q^k-1)/r}$
(see, for example, \cite{fr}, \cite{gal}).
If $P \in E( \F_q )$ and $k > 1$ then, since $\F_{q^k}$ is the extension
of $\F_q$  of minimal degree
which contains non-trivial $r$th roots of unity,
it follows that
$e( P, P ) = 1$ for the Tate pairing as well as the Weil pairing.

If $r \mid \#E( \F_ q )$ then the eigenvalues of
$\pi$ on $E( \overline{\F}_q )[r]$ are 1 and $q$.
Hence there is a basis $\{ P, Q \}$ for $E[r]$ such that
$\pi(P) = P$ and $\pi( Q ) = qQ$.
Now, $\{ P, Q \}$ forms a basis for the $r$-torsion
and so, by the same arguments as part 2 of Lemma~\ref{lemma-1},
$e( P, Q ) \ne 1$ for both Weil and Tate pairings.

Boneh observed (see \cite{friendly}, \cite{gal})
that the eigenspace $\langle Q \rangle$
of points with eigenvalue $q$ is equal to the set of
all points $R \in E( \F_{q^k} )[r]$ such that
$\Tr( R ) = 0_E$.
Boneh has also shown that $e( Q, Q ) = 1$
for the Tate pairing as well as the Weil pairing (see \cite{gal}).
We call $\langle Q \rangle$
the {\bf trace zero subgroup} and denote it
by $\T$.

\begin{lemma}
Let $E$ be an elliptic curve over $\F_q$. Let $r$ be
a large prime such that $r \mid \#E( \F_q)$ and $r \mid (q^k-1)$.
Define the basis $\{ P, Q \}$ as the eigenbasis for
Frobenius as above.
Let $S = aP + bQ \in E( \F_{q^k} )$ with $ab \ne 0$
and let $G = \langle S \rangle$.
Then the DDH problem in $G$ can be solved in polynomial time.
\end{lemma}

\begin{proof}
Consider $(S, uS, vS, wS)$.  Since $\Tr( S ) = kaP \ne 0_E$
and $b \ne 0$ we have $e( S, \Tr(S)) \ne 1$.
Hence, the DDH tuple $(S, uS, vS, wS)$
gives rise to the co-DDH tuple
\[
   ( S, uS, \Tr( vS ) = v \Tr(S), \Tr( wS) = w \Tr( S ))
\]
and, as we have seen, all co-DDH problems can be solved
using the Weil pairing.
\end{proof}

Hence, only two potentially hard DDH problems remain,
namely the subgroup $\langle P \rangle$ which is
the set of $r$-torsion points which are defined over
the field $\F_q$ and the trace zero subgroup $\T
\subset E( \F_{q^k} )[r]$.
Equivalently, these are the two eigenspaces in $E( \F_{q^k} )[r]$ for
the $q$-power Frobenius map.
In the ordinary case these problems seem to be hard.
For the remainder of the paper we consider the
supersingular case.

\section{Review of quaternion algebras}

We devote
this section to fixing the notation and briefly reviewing the
theory of quaternion algebras that we need in the sequel.

A quaternion algebra over a field $K$ is a central simple algebra
of rank $4$ over $K$. A quaternion algebra $B$ is division if
$B\not \simeq \M _2(K)$, or equivalently if $B^*=B\setminus \{ 0\}
$. If char$(K)\ne 2$, every quaternion algebra is of the form
$$
B=\left(\frac{m, n}{K}\right):=K+Ki+Kj+Ki j, \quad i^2=m, j^2=n, ij=-ji
$$
for some $m$, $n\in K^*$. The conjugation map on $B$ is $\overline
{a+bi+cj+dij} = a-bi-cj-dij$ and the reduced trace and norm on $B$
are $\Tr (\alpha ) = \alpha + \overline {\alpha }$ and $\n (\alpha
) = \alpha \cdot \overline {\alpha }$ for any $\alpha \in B$,
respectively.

We next present two different but equivalent versions of the
Skolem-Noether Theorem (see \cite{vig}).

\begin{proposition}
Let $B$ be a quaternion algebra over a field $K$.

\begin{enumerate}

\item Let $\sigma :B\rightarrow B$ be an automorphism of $B$ over
$K$. Then $\sigma (\alpha ) = \gamma ^{-1} \alpha \gamma $ for
some $\gamma \in B^*$.

\item Let $L/K$ be a quadratic field extension of $K$. Let $\phi ,
\psi :L\hookrightarrow B$ be two different immersions of $L$ into
$B$ over $K$. Then there exists $\gamma \in B^*$ such that
$\phi (\alpha )= \gamma ^{-1} \psi (\alpha )\gamma $ for all
$\alpha \in L$.

\end{enumerate}
\end{proposition}

Let $R$ be a Dedekind ring and let $K$ be its field of fractions.
Let $B$ be a quaternion algebra over $K$. We say that a place
$v\le \infty $ of $K$ ramifies in $B$ if $B\otimes K_v$ is a
division algebra over the completion $K_v$ of $K$ at $v$. A
classical theorem (see \cite{AlBa}, \cite{vig}, p.\,74) states
that there is a finite and even number of places of $K$ that
ramify in $B$.  Conversely, for any finite set $\{ v_1, ...,
v_{2r }\}$ of places
of $K$ of even cardinality, there exists a unique quaternion algebra
up to isomorphism
which ramifies exactly at the places $v_i$.

The reduced discriminant of $B$ is defined to be the product $D_B
= \prod \wp $ of all finite prime ideals of $R$ ramifying in $B$.

An element $\alpha $ in $B$ is integral over $R$ if $\Tr(\alpha
)$, $\n (\alpha )\in R$. Unlike number fields, the set of integral
elements in $B$ is not a subring of $B$ (for an example see
page 20 of \cite{vig}).

An order $\cR $ in $B$ over $R$ is a subring of $B$ of rank $4$
over $R$. We say that $\cR $ is maximal if it is not properly
contained in any other order of $B$. A left projective ideal $I$
of a maximal order $\cR $ is a locally principal sub-$\cR $-module
of $B$ of rank $4$ over $R$. Two projective left ideals $I$, $J$
of $\cR $ are linearly equivalent if $I = J\cdot \alpha $ for some
$\alpha \in B^*$. We let $\Pic _{R}(\cR )$ denote the set of
linear equivalence classes of left projective ideals of $\cR $
over $R$. The set $\Pic _{R}(\cR )$ is finite and its cardinality
$h_R(B)=\# \Pic _{R}(\cR )$ is independent of the choice of $\cR $.

The conjugation class of an order $\cR $ over $R$ is the set of
orders $[\cR ]=\{ \gamma ^{-1}\cR \gamma : \gamma \in B^*\}$,
which has infinite cardinality. There is however a finite number
$t_R(B)$ of conjugation classes of maximal orders in $B$ over $R$.

\begin{proposition}\label{type}

Let $K$ be the field of fractions of a Dedekind ring $R$ and let
$B$ be a quaternion algebra over $R$. Then

\begin{enumerate}
\item $h_R(B)\ge t_R(B)$.

\item If $K$ is a local field, then $h_R(B)=t_R(B)=1$.

\item If $K$ is a number field and $\mathfrak{M}$ is any ideal of $K$,
there exists an integral ideal $\mathfrak{N}$ of $R$, $(\mathfrak{M},
\mathfrak{N}) = 1$, such that
$h_{R[\frac{1}{\mathfrak{N}}]}(B)=t_{R[\frac{1}{\mathfrak{N}}]}(B)=1$.

\end{enumerate}
\end{proposition}

{\em Proof. } The first two statements can be found in \cite{vig},
p.\,26 and Ch.\,II respectively. As for the third, let $\cR $ be a
maximal order of $B$ and let $\{ I_1$, ..., $I_{h_{R}(B)}\} $ be a
full set representatives of projective left ideals in $\Pic_R(\cR
)$. It follows from \cite{Ro}, p.\,5 that $I_i$ can be chosen such
that $\mathfrak{N} = \n (I_1)\cdot ... \cdot \n (I_{h_{R}(B)})$ is
coprime to $\mathfrak{M}$. Since $I_i$ are invertible in $\cR
[\frac{1}{\mathfrak{N}}]$, we have that
$h_{R[\frac{1}{\mathfrak{N}}]}(B)=1$. By $(1)$ we also have
$t_{R[\frac{1}{\mathfrak{N}}]}(B)=1$. $\Box $

Let $B=(\frac{m, n}{K}):=K+Ki+Kj+Ki j, i^2=m, j^2=n, ij=-ji$ and
let $\cR $ be a maximal order in $B$ over $R$. Two questions that
naturally arise in several contexts and that we encounter in the
proof of Theorem \ref{main} are the following:
\begin{enumerate}
\item Do there exist elements $\pi $, $\psi \in \cR $ such that
$\pi ^2=m, \psi ^2=n, \pi \psi =-\psi \pi $?

\item Fix $\pi \in \cR $ such that $\pi ^2=m$ (if there is any).
Does there exist $\psi \in \cR $ such that $\psi ^2=n, \pi \psi
=-\psi \pi $?
\end{enumerate}

These questions were considered in the appendix to \cite{Ro1}. We
state here a partial answer which will suffice for our purposes.

\begin{proposition}\label{basis}
Let notation be as above.
\begin{enumerate}
\item If $t_R(B)=1$, then there exist $\pi $, $\psi \in \cR $ such
that $\pi ^2=m, \psi ^2=n, \pi \psi =-\psi \pi $.

\item Fix $\pi \in \cR $ such that $\pi ^2=m$. If $t_R(B)=1$ and
$\OO = R[\sqrt{m}]\subset K(\sqrt{m})$ is locally a discrete
valuation ring at the places $v\nmid D_B$ of class number $h(\OO)
= 1$, then there exists $\psi \in \cR $ such that $\psi ^2=n, \pi
\psi =-\psi \pi $.

\end{enumerate}

\end{proposition}

{\em Proof. } Part $(1)$ follows from \cite{Ro1}, Proposition 5.1.
As for $(2)$, let $\mathcal E(m)$ denote the set of embeddings
$i:R[\sqrt{m}]\hookrightarrow \cR $ over $R$ up to conjugation by
elements in the normalizer group $\mathrm{Norm}_{B^*}(\cR )$.
Since $\pi \in \cR $, $\mathcal E(m)$ is non empty. Eichler proved
that $\mathcal E(m)$ is a finite set. More precisely, we have from
our hypothesis and \cite{vig}, Theorem 3.1 on p.\,43 and Theorem
5.11 on p.\,92, that in fact $\# \mathcal E(m)=1$. It now follows
from \cite{Ro1}, Proposition 5.7 and its remark below that there
exists $\psi \in \cR $ such that $\psi ^2=n, \pi \psi =-\psi \pi
$. $\Box $

Let $B$ be a quaternion algebra over $\Q $. We say that $B$ is
definite if $\infty $ ramifies in $B$, that is, if $B\otimes \R
=\mathbb H$ is the algebra of real Hamilton quaternions.
Equivalently, $B$ is definite if and only if $D_B$ is the product
of an odd number of primes. Otherwise $B\otimes \R =\M _2(\R )$
and we say that $B$ is indefinite.

If $B$ is indefinite then $h_{\Z }(B)=t_{\Z }(B)=1$. Otherwise,
$h_{\Z }(B)$ and $t_{\Z }(B)$ can explicitly be computed as in
\cite{vig}, p.\,152. When $D_B=p$ is prime, the class number
$h_{\Z }(B)$ is the number of isomorphism classes of supersingular
elliptic curves over $\bar \F_p$ and $t_{\Z }(B)$ is the number of
isomorphism classes of supersingular elliptic curves up to $\Gal
(\bar \F_p/\F_p)$-conjugation.

Let $\Q _v$ be a local completion of $\Q $ at a place $v\le \infty
$. The Hilbert symbol over $\Q _v$ is a symmetric bilinear pairing
$$
(\quad , \quad )_v:\Q _v^*/\Q _v^{*2}\times \Q _v^*/\Q
_v^{*2}\rightarrow \{ \pm 1 \}
$$
which may be defined as $(m, n)_v = 1$ if the quaternion algebra
$(\frac{m, n}{\Q _v})\simeq \M _2(\Q _v)$ and $(m, n)_v = -1$
otherwise.

In practice, the Hilbert symbol is computed as follows. For
$v=\infty $, $(m, n)_{\infty } = -1$ if and only if $m<0$ and
$n<0$. For any odd prime $p$, $(m, n)_p$ can be computed by using
the multiplicative bilinearity of the pairing and the following
three properties:
\begin{itemize}
\item $(-p, p)_p=1$ \item $(m, n)_p =1$ if $p\nmid 2 m n$ \item
$(m, p)_p = (\frac {m}{p})$ is the Legendre quadratic symbol for
any $p\nmid m$.
\end{itemize}

Finally, the Hilbert symbol at $2$ follows from the equality $\prod _v
(m, n)_v = 1$.

\section{Supersingular curves and distortion maps}
\label{main-theorem}

In the next sections we restrict attention to supersingular
curves. As is known (see, for example, \cite{silv} Theorem V.3.1 and
\cite{gross}), an elliptic curve $E$ over a finite field $\F _q$
is supersingular if and only if $\End_{\overline{\F_q}}(E)\otimes
\Q $ is a quaternion algebra over $\Q $ of reduced discriminant
$p$.

Verheul \cite{verheul} was the first to propose using non-rational
endomorphisms to solve DDH problems. Let $P \in E( \F_{q^k} )$ be
a point of order $r$. If $\psi \in \End( E )$ is such that
$\psi(P) \not\in \langle P \rangle$ then $\{ P, \psi(P) \}$ is a
generating set for $E[r]$ and so $e_r( P, \psi(P)) \ne 1$.  It
follows that DDH problems in $\langle  P \rangle$ can be solved.
Verheul called such endomorphisms {\em distortion maps}.

Originally, distortion maps were exclusively used to map points
defined over $\F_q $ to points defined over $\F_{q^k}$. In other
words, the focus had been on the 1-eigenspace for the Frobenius
map on $E[r]$. Theorem 5 of Verheul \cite{ver-full}
states that a suitable distortion map exists for every point
$P \in E[r]$ when $E$ is supersingular. The proof of Theorem 5 of
\cite{ver-full} is not constructive, and it seems difficult to
obtain an algorithm for finding a distortion map using that
approach.

In Theorem \ref{main} below we obtain
an analogous result to that in \cite{ver-full} using completely
different techniques.  We can then give in Section
\ref{algorithm} an algorithm for constructing a distortion map for any
supersingular curve.



\begin{lemma}
\label{lemma-deg}
Let $E$ be a supersingular elliptic curve over $\F_q$
and let $\psi$ be an endomorphism.  Let $P$
be an element of one of the eigenspaces of the
$q$-power Frobenius map $\pi$.  Then
$\psi$ maps $P$ outside $\langle P \rangle$
if and only if
\[
   P \not\in \ker ( \psi \pi - \pi \psi ) .
\]
\end{lemma}

\begin{proof}
Suppose $\pi( P ) = [m]P$ for
some $m$ (indeed, either $m=1$ or $m=q$).
Now, $\psi(P) $ also in the eigenspace means
$\pi \psi(P) = [m] \psi( P ) = \psi ( [m]P) = \psi \pi( P) $.
In other words, $P \in \ker( \psi \pi - \pi \psi )$.
The converse is similar.
\end{proof}

\begin{theorem}\label{main}
Let $E$ be a supersingular curve over $\F_q$, $q=p^a$.
Suppose $k > 1$ and let $r\mid  \# E(\F_q)$,
$r\not = p$, $r>3$, be a prime.
Let $\pi$ be the $q$-power Frobenius map and
let $P \in E(\F_{q^k})$ be in a $\pi$-eigenspace.
Then there exists a
distortion map $\psi$ on $E$ which maps $P$ outside
$\langle P \rangle$.
\end{theorem}

\begin{proof}
By Lemma~\ref{lemma-deg}, to prove the result
it is enough to prove that
there exists $\psi \in \End (E)$ such
that $r \nmid \deg (\pi \psi - \psi \pi )$.

Let $P(T)=T^2 - t T + q$ be the characteristic polynomial of the
$q$-power Frobenius element $\pi $ acting on $E$. Since $k>1$, we
know (see for example \cite{water} or \cite{gal}, Theorem I.20) that
$P(T)$ is irreducible and so its roots
generate a quadratic field of $\Q $.

The endomorphism ring $\cR = \End (E)$ is a maximal order in the
quaternion algebra $B_p=\End (E)\otimes \Q $, which ramifies
exactly at $p$ and $\infty $ \cite{gross}. The ring $\End _{\F
_q}(E)$ is an order in the quadratic field $\Q (\pi )=\End
_{\F_q}(E)\otimes \Q \simeq \Q (\sqrt{t^2-4q})$, naturally
embedded in $\cR $. Let $\pi _0=2\pi -t\in \Q (\pi )$, which
satisfies $\Tr (\pi_0)=0$ and $\n(\pi_0)= - \pi_0^2 = 4q-t^2$.

There is a morphism of $\Q $-vector spaces
$$
\begin{matrix}
c_{\pi }: & B_p & \rightarrow & B_p \\
  & \psi &\mapsto & \pi \psi -\psi \pi
\end{matrix}
$$
with $\ker (c_{\pi }) = \Q (\pi )$.

Let $s\in \Z $. We remark that there exists an element $\psi _0\in
B_p$ such that $\psi _0^2=-s$ and $\pi _0 \psi _0 = -\psi _0 \pi
_0$ if and only if
$$
B_p\simeq \left(\frac{t^2-4q,-s}{\Q }\right) \qquad \qquad (\dag ).
$$

Indeed, one direction is immediate. The other
implication follows from
the Skolem-Noether Theorem: If $B_p=\Q +\Q i+\Q j+\Q ij =
(\frac{t^2-4q,-s}{\Q })$, there exists $\gamma \in B_p^*$ with
$\pi_0=\gamma ^{-1} i \gamma $ and we may take $\psi_0=\gamma
^{-1} j \gamma $.

Note that since the discriminant of $B_p$ is $p$, condition
$(\dag) $ for a given $s$ can be checked by computing a finite number of
local Hilbert symbols. Moreover, since $B_p\otimes \R $ is a
division algebra, necessarily $s>0$.

Let $s\in \Z $ be such that $(\dag)$ holds. By Proposition \ref{type}, $(3)$
there exists an integer $N_0$ coprime to $r$ such that $t_{\Z
[\frac{1}{N_0}]}(B_p)=1$.
Similarly, there exists an integer $N_1$ coprime to $r$ such that
$\Z [\frac{1}{N_1}, \sqrt{t^2-4q}]$ is locally a discrete
valuation ring at all primes $\ell \not = p$ and $h(\Z
[\frac{1}{N_1}, \sqrt{t^2-4q}])=1$.   Indeed, this is accomplished
by considering a system of representatives $J_1, ..., J_{h(\Q
(\sqrt{t^2-4 q}))}$ of classes of ideals in the quadratic field
$\Q (\sqrt{t^2-4 q})$ such that $r\nmid \mathrm{N}_{\Q
(\sqrt{t^2-4 q})/\Q }(J_i)$ and taking $N_1 = 2\cdot \prod
\mathrm{N}_{\Q (\sqrt{t^2-4 q})/\Q }(J_i)$.

By Proposition \ref{basis}, there exists $\psi_0\in \cR
[\frac{1}{N_0\cdot N_1}]$ such that $\psi _0^2=-s$ and $\pi _0
\psi _0 = -\psi _0 \pi _0$. Hence $N \psi_0\in \cR $ for some
integer $N$ supported at the primes dividing $N_0\cdot N_1$ and
thus coprime to $r$.
The endomorphism $N \psi_0$ will be the distortion map
we are looking for (this is all assuming that $(\dag)$ holds).

Since $\pi \psi _0 -\psi _0 \pi = \pi_0 \psi _0$ and $\pi (\pi_0
\psi _0) -(\pi _0 \psi _0) \pi = (\frac{\pi_0-t}{2}) (\pi_0 \psi
_0) -(\pi _0 \psi _0) (\frac{\pi_0-t}{2}) =(t^2-4 q)\psi _0$, it
readily turns out that $\mathrm{Im}(c_{\pi }) = \Q \cdot \psi _0 +
\Q \cdot \pi _0 \psi_0$ and $c_{\pi }(\cR )\supseteq c_{\pi }(\Z +
\Z \pi_0 + N \Z \psi_0 + N \Z \pi \psi_0) = (t^2-4q) N \Z \psi_0 +
N \Z \pi_0 \psi_0$. Moreover, the degree of the isogenies
$(t^2-4q) N \psi_0$ and $N \pi_0 \psi_0$ on $E$ are computed in terms
of the reduced norm in the quaternion algebra $B_p$ as $\deg
((t^2-4q) N \psi_0)=(t^2-4 q)^2 N^2 \n (\psi _0)=(t^2-4 q)^2 N^2
s$ and $\deg (N \pi_0 \psi_0) = N^2 \n(\pi_0) \n (\psi
_0)=(4q-t^2) N^2 s$.
Hence,
\[
   \deg( \pi (N \psi_0) - (N \psi_0) \pi ) = N^2 (4q - t^2) s
\]
is coprime to $r$ as desired.

It remains to give choices of $s$ for which condition
$(\dag)$ is satisfied.
According to a theorem of Waterhouse \cite{water}
the possible values of the trace of the
Frobenius endomorphism are $t=0, \pm p^{a/2}, \pm 2 p^{a/2}$ and
$\pm p^{(a+1)/2}$. Recall that we can exclude the value $t=\pm
2p^{a/2}$ because we are assuming $k>1$. Hence, the only possible
prime factors of $4 q-t^2$ are $2$, $3$ and $p$, and in order to
prove the claim, it suffices to show that $B_p\simeq
(\frac{t^2-4q, -s}{\Q })$ for either $s=1$ or for some prime $s$,
$s\not = r$.

The following table lists,
for each of the possible values of $t$, a choice of $s$
such that condition $(\dag )$ holds:

\begin{center}
\begin{tabular}{|c|c|}
  \hline
    If $t=0$, $a$ is odd, $p\not \equiv 1 \mod 4$ & $s=1$ \\
  \hline
  If $t=0$, $a$ is odd, $p\equiv 1 \mod 4$ & Any prime $s\equiv 3\mod 4$  and split in $\Q (\sqrt{-p})$ \\
  \hline
  If $t=0$, $a$ is even & $s=p$ \\
\hline
  If $t=\pm p^{(a+1)/2}$ & $s=1$\\
\hline
  If $t=\pm p^{a/2}$ & $s=p$ \\
  \hline
\end{tabular}
\end{center}

This table is checked by computing relevant Hilbert symbols.
We give details of the argument for the first two rows of
the table.
Assume that $t=0$ and $a$ is odd. We have that $(-4p^a,
-s)_{\ell } = (-p, -s)_{\ell }$ for all primes $\ell $ and $(-p,
-s)_{\ell }=1$ for all finite primes $\ell \nmid 2 p\cdot s$.
Moreover, we have $(-p, -s)_{\infty } = -1$ if  and only if $s>0$.

If $p\not \equiv 1 \mod{4}$, $p\ne 2$, then $(-p, -1)_{p} = (p,
-1)_{p} = (\frac{-1}{p}) = -1$. Since we have that $p$ and $\infty
$ ramify in $(\frac{-p, -1}{\Q })$ and the number of ramifying
places must be even, we have that $(-p, -1)_{2}=1$. Hence
$(\frac{-p, -1}{\Q })$ is the quaternion algebra of discriminant
$p$ and $B_p\simeq (\frac{-p, -1}{\Q })$.

Similarly, if $p=2$, it holds that $B_2\simeq (\frac{-2, -1}{\Q
})$.

If $p\equiv 1 \mod 4$ and $s$ is a prime $s\equiv 3\mod{4}$
and split in $\Q (\sqrt{-p})$ (i.e.,
$(\frac{-p}{s}) = 1$), then $(-p, -s)_{p} = (p,
-s)_{p} = (\frac{-s}{p}) = -1$ and $(-p, -s)_{s} = (-p, s)_{s} =
(\frac{-p}{s}) = 1$. Hence $B_p\simeq (\frac{-p, -s}{\Q })$.

Note that
the Theorem of \v{C}ebotarev implies there are infinitely many
suitable primes $s$ for line two of the table, hence we
can always choose one which is not divisible by $r$.

We leave the checking of the remaining cases of the table above to
the reader; remember that line three of the table only applies
to $p=2$ or $p \equiv 3 \mod{4}$, that line four of the
table only applies to $p=2, 3$ and that line five of
the table only applies when $p=3$ or $p \equiv 2 \mod{3}$.

This completes the proof.
\end{proof}

\begin{remark}
It follows from the above proof that Theorem \ref{main} is also
valid for $r=3$ unless $p=3$ or $t=\pm p^{a/2}$. The statement is
valid for $r=2$ precisely when $p\ne 2$ and $t=\pm p^{a/2}$ or
when $p=3$ and $t=\pm p^{(a+1)/2}$.
\end{remark}

\section{An algorithm for constructing distortion maps}
\label{algorithm}

The aim of this section is to derive from the proof of Theorem
\ref{main} an algorithm for constructing a distortion map on a
supersingular curve over a field of characteristic $p$.

One might expect the first step of such an algorithm to involve
computing a basis for the endomorphism ring using Kohel's
algorithm \cite{kohel} (which runs in exponential time). In fact,
we argue that this is not required. Instead we reflect upon how
one would obtain a usable supersingular elliptic curve.
It is
known that for all finite fields $\F_q$ there is a supersingular
elliptic curve $E$ defined over $\F_q$ (and in general, there will
be many non-$\F_q$-isomorphic such curves).  We
claim that all the curves which could potentially
be used in practice arise as reductions
of CM curves in characterstic zero of small class number.


To justify our claim, consider the
following three candidate methods to find a supersingular
curve over a finite field.
\begin{enumerate}
\item Using the complex multiplication (CM) method.

\item Constructing curves over fields of small characteristic.
For example $y^2 + y = f(x)$ over $\F_{2^m}$ is always supersingular.

\item Choosing random curves over $\F_p$ or $\F_{p^2}$
and counting points until a supersingular curve is found.
\end{enumerate}

The third method method is not useful as the probability of
success is negligible.  The number of isomorphism classes of supersingular
curves over $\F_p$ is equal to $h_{-4p} + h_{-p}$
(where $h_{D}$ is the class number of the order of discriminant
$D$, and $h_{-p} = 0$ if $p \equiv 1 \mod{4}$, for details
see Gross~\cite{gross}).
By the Brauer-Siegel theorem (more details below)
this number is roughly $p^{1/2}$ and so the
probability of a randomly chosen elliptic curve over $\F_p$
being supersingular is negligible.
Similarly, the number of isomorphism classes of supersingular
curves over $\F_{p^2}$ is $p/12$
(see Theorem V.4.1 (c) of \cite{silv}) and so the probability
of a random curve over $\F_{p^2}$ being supersingular is negligible.

The second method restricts attention to a
very small number of isomorphism classes (and hence $j$-invariants).
In the example given, the curves all have $j=0$.
Hence, all these curves can be treated as twists of reductions
of curves in characteristic zero, and these curves can be
chosen to be CM curves.
Hence the second method is essentially a special case of the CM method.

The CM method works in the following setting.
Let $E$ be an elliptic curve over a number field
$F$ with complex multiplication by an order $\OO$ in an imaginary
quadratic field $K=\Q( \sqrt{-d} )$. Let $p$ be a rational prime
which does not split in $\OO$ and let $\pp$ be a prime
of $F$ above $p$.  Then by the Deuring reduction theorem,
$\tilde{E} = E \mod{\pp}$ is a supersingular elliptic
curve over the residue field $k$ of $F$ at $\pp $.
The main step of the CM method  is to construct the
ring class polynomial of the order $\OO$ (which has
degree $h_\OO$, the class number of the order)
and to find a root of it in characteristic $p$.
This process has exponential complexity in the
class number $h_\OO$ and can only be applied in
practice when $h_\OO$ is relatively small.

It would be very interesting to have an alternative construction
for supersingular curves.  This open problem is also raised in
Section 4.1 of Verheul \cite{ver-full}.

\begin{proposition}\label{red}

Let $E/F$ be an elliptic curve defined over a number field $F$
with complex multiplication by an order $\OO$ of discriminant
$D$ in an imaginary quadratic field $K=\Q(\sqrt{D})$.
Assume that $K\not\subset F$.
Let $p$ be a prime for which $E$ has good and supersingular reduction.
Let $\pp $ be a prime ideal of $F$ above $p$. Let $\tilde E$
over $k=\F_{p^m}$ be the reduction mod $\pp $ of $E$.
Let $\pi$ be the $p^m$-Frobenius map on $\tilde{E}$.
Suppose $r | \# \tilde{E}( \F_{p^m} )$ is a prime
such that $r > 3$ and $r \nmid pD$.

Let $d > 0$ be such that $\sqrt{-d} \in \OO$.
Let $\Psi \in \End (E)$ satisfy  $\Psi^2=-d$.
Let $\psi \in \End_{\bar \F_p}(\tilde E)$ be the reduction mod $\pp$ of $\Psi$.
Then $\psi$ is a suitable distortion map for points $P \in \tilde{E}[r]$
which lie in a $\pi$-eigenspace.
\end{proposition}

\begin{proof} Note that since $K\not \subset
F$, $H=F\cdot K$ is a quadratic extension over $F$. We know by the
theory of complex multiplication that the minimal field of
definition of the endomorphisms of $E$ is $H$ and it follows that,
if we let $\sigma \in \Gal (H/F)$ be a non trivial element, then
$\Psi^{\sigma } = -\Psi$. Let $\tilde k$ be the residue field of a prime
ideal in $H$ above $\pp $. The natural Galois action of $\Gal
(H/F)$ on $\End _{H}(E)\otimes \Q $ descends to an action of $\Gal
(\tilde k/k)$ on $\End _{\tilde k}(E)\otimes \Q \simeq B_p$. If we
let $\tilde \sigma $ denote a generator of $\Gal (\tilde k/k)$, we
have that $\psi^{\tilde{\sigma }}=-\psi$ due to the compatibility of the
Galois action.

The Galois automorphism $\tilde{\sigma }$ acts on the quaternion
algebra $B_p$ as an automorphism $\tilde{\sigma }:B_p\rightarrow
B_p$. By the Skolem-Noether Theorem, $\alpha ^{\tilde{\sigma }} =
\gamma \alpha \gamma ^{-1}$ for some $\gamma \in B_p^*$, which is
uniquely determined as an element of $B_p^*/\Q ^*$. Since $\pi
^{\tilde{\sigma }} = \gamma \pi \gamma ^{-1} = \pi $ because $\pi
\in \End_{k}(E)$, we deduce that $\gamma \pi = \pi \gamma $ and
hence $\gamma \in \Q (\pi )$. Since $\psi^{\tilde{\sigma }}=\gamma \psi
\gamma ^{-1}=-\psi$, it follows that $\Tr(\gamma \psi)=\gamma \psi +
\overline{\gamma \psi} = -\psi \gamma  + \overline{\psi}\overline{\gamma }
= \gamma \psi - \psi \overline{\gamma } = -\Tr(\gamma ) \psi \in \Z $. Hence
$\Tr(\gamma )=0$ and $\gamma = \pi $ in $B_p^*/\Q^*$. Thus $\pi \psi
= -\psi \pi $ and so $\psi \pi - \pi \psi = 2 \psi \pi$ is an
isogeny of degree $4(4p - t^2) d$.

Now let $P \in \tilde{E}[r]$ be in a $\pi$-eigenspace.
We apply arguments used in the proof of Theorem~\ref{main}.
Since $r >3$ and $r \nmid pd$ we have
that $P \not\in \ker( \psi\pi - \pi \psi)$.
Therefore $\psi(P)$ is independent from $P$.
\end{proof}

We can now present our algorithm.  The input is a supersingular
elliptic curve $\tilde{E}$ over a finite field
$\F_q$ where $q=p^m$.
We also assume that an order $\OO \subset \End(\tilde{E})$ of
class number $h_\OO$ is specified.
Note that by the Brauer-Siegel theorem
(see Theorem XVI.5 of Lang~\cite{lang-ant}, for non-maximal
orders also see Theorem 8.7 of \cite{lang})
we have that the discriminant $D_\OO$ of $\OO$ is
$O( h_\OO^{2+\epsilon})$.
The notation $D_\OO = O( h_\OO^{2 + \epsilon})$ means
that for every $\epsilon > 0$ there is a constant
$c_\epsilon$, which depends on $\epsilon$,
such that $D_\OO \le c_\epsilon h_\OO^{2 + \epsilon}$
for all $\OO$.

\vskip 0.4cm

\noindent
{\bf ALGORITHM 1: Construction of a distortion map on $\pmb{\tilde E}$:}

\begin{enumerate}
\item  Let $\OO$ be an order in $\End( \tilde{E} )$ of class
number $h_\OO$.
Compute the discriminant $D$ of $\OO$.  Hence compute an integer
$d>0$ of size $O(D)$ such that $\sqrt{-d} \in \OO$ (for example, we can take
$d=-D$). Denote $\sqrt{-d}$ by $\psi$, so that $\psi$ is a
$d$-isogeny.

\item Factor $d$ as $\prod_{i=1}^n l_i$
(where $l_i$ are not necessarily distinct primes). Then $\psi$
is a composition $\psi_1 \cdots \psi_n$ of prime degree
isogenies (and each $\psi_i$ will be defined over
$\F_{q^2}$).

\item Use Galbraith's algorithm \cite{gal-isog} to construct a
tree of prime degree isogenies between $j$-invariants of
supersingular elliptic curves in characteristic $p$.
The tree starts with vertex $j(\tilde{E})$ and the process
terminates when this vertex is revisited by a non-trivial isogeny.
Since we know there is a non-trivial isogeny $\psi$ of degree
$d$ we should select only the primes $l_i$ as found in step (2).

\item Construct the isogeny $\psi$ on $\tilde{E}$ explicitly as the
composition of isogenies $\psi_i$.
Each isogeny $\psi_i$ can be computed from the
$j$-invariants of the corresponding elliptic curves
using methods of Elkies \cite{elkies} and V\'elu
\cite{velu}. Usually it is also necessary to construct an additional
isomorphism between the image of the final isogeny $\psi_n$
and the elliptic curve $\tilde{E}$.
All these calculations will be performed over $\F_{q^2}$.
\end{enumerate}

\vskip 0.4cm

By Proposition \ref{red}, the endomorphism $\psi$ will be
a suitable distortion map.  Hence the algorithm is clearly
correct.

We now roughly analyse the complexity of the algorithm.
We assume a unit
cost for operations in the field of definition $\F_{q}$ of
$\tilde{E}$.  We express the complexity in terms of the
class number $h = h_\OO$.  For further details of the
complexity analysis of algorithms like this see Elkies \cite{elkies} and
Galbraith~\cite{gal-isog}.

\begin{enumerate}
\item Step one is essentially trivial.  Since $D$ is $O( h^{2 + \epsilon})$
the complexity of this step is $O( h^{2 + \epsilon})$.

\item Factorisation can be easily done in time $O( \sqrt{d})$
which is $O( h^{1 + \epsilon})$.  The number $n$ is $O( \log(h) )$
while the primes themselves are $O( d )=O( h^{2 + \epsilon})$.

\item There are $n =O( \log(h))$ iterations of  the process.
Each step requires computing the $l$-th modular polynomial
$\Phi_l(x,y)$ (which has degree $l+1$ in each
variable and takes $O( l^3)$ operations to compute)
and finding the roots of $\Phi_l( j, y)$ in $\F_{q^2}$
(which takes $O( l  \log(q))$ operations).
The total cost of this stage in the worst
case is therefore,
$O( \log(h) (h^{6+\epsilon} + h^{2+\epsilon} \log(q)))$.
The space requirement for the tree is $O( \log(h))$.

\item Finding the path in the tree takes time $O( \log(h))$.
For each $l$-isogeny in the composition,
Elkies' algorithm requires $O( l^3 )$ operations and
V\'elu requires $O(l)$ operations.
Computing the isomorphism is trivial.
Hence the total cost of explicitly computing the isogeny $\psi$
is $O( \log(h) h^{6 + \epsilon} )$ operations.
\end{enumerate}

To conclude, it is clear that step 3 is the dominant step.
The total complexity of the algorithm is
$O( \log(h) (h^{6+\epsilon} + h^{2+\epsilon} \log(q)))$.
Since we can only construct curves for which $h$ is
bounded by a polynomial function, this is therefore a
polynomial time algorithm on families of curves which have
been constructed in any practical setting.



%
%

\section{Standard examples}
\label{sec-7}

In the previous sections we showed the existence of non-rational
endomorphisms $\psi$ with a certain property (namely, that $\psi(
\pi( Q )) \ne \pi (\psi(Q))$ for points of order $r$ which are
in a Frobenius eigenspace).
In practice there are a small number of examples of supersingular curves
which are widely used, and popular distortion maps are already
known in these cases.
In this section we recall these familiar examples
and show that they satisfy the above property.

Table~\ref{table1} gives the list of curves studied.
These curves have been considered by several
authors (for example, Verheul~\cite{verheul} and Galbraith~\cite{galss}).
Note that in all cases we have
$j(E) = 0$ or 1728.
This table does not list all possible variations
of distortion maps.  For instance, Barreto has
suggested using
\[
   \psi( x,y) = (x + \zeta_3^2, y + \zeta_3 x + t)
\]
where $ t^2 + t = \zeta_3$
in the case of characteristic 2 and $k=4$.

\begin{table}
\begin{tabular}{|l|l|} \hline
 $k$ & Elliptic curve data \\
 \hline
 2 & $E : y^2 = x^3 + a$ over $\F_{p}$ where $p \equiv 2 \mod{3}$, $p>2$ \\
 & $\#E( \F_p ) = p+1$ \\
 & Distortion map $(x,y) \mapsto (\zeta_3 x, y )$
where $\zeta_3^3 = 1$.\\
\hline
 2 &  $y^2 = x^3 + a x$ over $\F_{p}$ where $p \equiv 3 \mod{4}$ \\
 & $\# E( \F_p) = p+1$. \\
 & Distortion map $(x,y) \mapsto (-x, i y)$ where $i^2 = - 1$. \\
\hline
 3 & $E : y^2 = x^3 + a$ over $\F_{p^{2}}$ where \\
  & $p \equiv 5 \mod{6}$ and $a \in \F_{p^2}$, $a \not\in \F_p$ is a square
    which is not a cube. \\
 & $\#E( \F_{p^2} ) = p^{2} - p + 1$. \\
 &  Distortion map $(x,y) \mapsto (\gamma^2 x^p ,
  b y^p /b^p )$  \\
  & where $a = b^2$ ($b \in \F_{p^2}$) and $\gamma \in \F_{p^6}$ satisfies $\gamma^3 = b/b^p$. \\
\hline
  4 & $y^2 + y = x^3 + x + b$ over $\F_{2}$ \\
 & Distortion map $(x,y) \mapsto (\zeta_3 x + s^2, y + \zeta_3 s x + s )$ \\
 & where $s \in \F_{2^4}$ satisfies $s^2 + \zeta_3 s + 1 = 0$. \\
\hline
 6 & $ y^2 = x^3 + a x + b$ over $\F_{3}$. \\
 & Distortion map
 $(x,y) \longmapsto ( \alpha - x , i y)$
 where $i \in \F_{3^2}$ and $\alpha \in \F_{3^3}$ \\
 & satisfy $i^2 = -1$ and $\alpha^3 + a \alpha - b = 0$. \\
\hline
\end{tabular}
\caption{Popular distortion maps.}
\label{table1}
\end{table}

%
%
%

\begin{proposition}
Let $E$ be a supersingular curve over $\F_q$
from Table~\ref{table1}
where $q$ is a power of $p > 3$.
Let $\pi$ be the $q$-power Frobenius.
Suppose $ r \mid \# E( \F_q )$
and $r > 3$.
Then the distortion map $\psi$ listed in the
table satisfies $r \nbar \deg( \pi \psi - \psi \pi )$.
\end{proposition}

\begin{proof}
Consider first the case when
$E$ is the curve $ y^2 = x^3 + ax$ over $\F_p$
with $k=2$ and
with the distortion map $\psi : (x,y) \mapsto (-x, iy)$.
Clearly, $\psi^2 = -1$ and this case is covered by
Proposition~\ref{red}. One can also give a direct proof.

Now consider the case
$E: y^2 = x^3 + a$ with $k=2$ over
$\F_p$ ($p \equiv 2 \mod{3}$) and with
the distortion map
$\psi( x,y) = (\zeta_3 x, y)$.
In this case we have $\psi^3 = 1$ and
so Proposition~\ref{red} does not apply.
A variant of Proposition~\ref{red} which handles this
case can be proved, but instead
we give the following direct argument.
Let $Q = (x,y) \in E[r]$.
Since $r > 3$ we have $x \ne 0$.
We have $\pi \psi (Q) = \pi( \zeta_3 x, y ) = ( \zeta_3^2 \pi(x), \pi(y) )$
while $\psi \pi(Q) = (\zeta_3 \pi(x), \pi(y)) \ne \pi \psi(Q)$.
Clearly, $Q \not\in \ker( \pi \psi - \psi \pi)$ and
the result follows.

Finally, consider the case $k=3$ with $E : y^2 = x^3 + a$.
Since $\gamma^2 \not\in \F_{p^2}$ we have $\pi(\gamma^2) \ne \gamma^2$.
The $x$-coordinate of $\psi \pi(Q)$ is
$\gamma^2 \pi(x)$ while the $x$-coordinate
of $\pi \psi(Q)$ is $\pi( \gamma^2 x ) = \pi(\gamma^2) \pi(x)$.
Since $r >3$ we have $x\ne 0$, and so the $x$-coordinates
are not equal.  The result follows.
\end{proof}

\begin{proposition}
Let $E$ be a supersingular curve over $\F_q$
from Table~\ref{table1}
where $q$ is a power of 2.
Let $\pi$ be the $q$-power Frobenius map.
Suppose $ r \mid \# E( \F_q )$ is such that $r > 1$.
Then the distortion map $\psi$ listed in the
table satisfies $r \nbar \deg( \pi \psi - \psi \pi )$.
\end{proposition}

\begin{proof}
The relevant curve is $E : y^2 + y = x^3 + x + b$
with distortion map
$\psi( x, y ) = (\zeta_3 x + s^2, y + \zeta_3 s x + s)$
where $\zeta_3^3  = 1$ and
$s^2 + \zeta_3 s  + 1 = 0$.
Since $\psi^3 = 1$ we cannot apply Proposition~\ref{red}
so we give a direct argument.

If $\pi \psi(Q) = \psi \pi(Q)$ then $\pi^2 \psi(Q) = \psi \pi^2(Q)$
so it is enough to prove that the latter equality does not hold.
Suppose $q = 2^m$ where $m$ is odd (otherwise $k < 4$).
Clearly, $\pi^2$ fixes $\F_{q^2}$ and so
$\pi^2( \zeta_3) = \zeta_3$.
Now $\pi^2$ does not fix $s \in \F_{q^4}$ so, by inspection
of the minimal polynomial, $\pi^2(s) = s + \zeta_3$.

Let $Q =(x,y)\in E[r]$. Then the $x$-coordinate of
$\pi^2 \psi( Q ) $ is
$\pi^2( \zeta_3 x + s^2 ) = \zeta_3 \pi(x^2) + s^2 + \zeta_3^2$
while the $x$-coordinate of $\psi \pi^2(Q)$ is
$\zeta_3 \pi^2(x) + s^2$.
The result follows.
\end{proof}

\begin{proposition}
Let $E$ be a supersingular curve over $\F_q$
from Table~\ref{table1}
where $q$ is a power of 3.
Let $\pi$ be the $q$-power Frobenius map.
Suppose $ r \mid \# E( \F_q )$
and $r > 1$.
Then the distortion map $\psi$ listed in the
table satisfies $r \nbar \deg( \pi \psi - \psi \pi )$.
\end{proposition}

\begin{proof}
Clearly, $\psi^2 = -1$ and Proposition~\ref{red} applies
(take $F$ to be a cubic extension of $\Q$).
There is also an easy direct proof.
%
%
\end{proof}

\section{Distortion maps which are not isomorphisms}
\label{not-iso}

By Theorem III.10.1 of \cite{silv}
there are non-trivial automorphisms only
when $j(E) = 0$ or 1728 (in particular, when
the endomorphism ring has a subring isomorphic to either
$\Z[ i ] $ or $\Z[ \zeta_3 ]$, both of which are rings with
non-trivial units).
Hence, we cannot expect distortion maps to be
automorphisms in all cases.

Even in the cases $j=0, 1728$ we see that the value
$s=1$ cannot always be taken in the proof of
Theorem~\ref{main}.  This indicates why
the $k=3$ example in characteristic $p$
(with $t=p^{a/2}$)
does not admit an automorphism.

The aim of this section is to give some examples
of these distortion maps.
For the first example we use Algorithm 1.
For the second example we use an ad-hoc technique
which shows that Algorithm 1 is not optimal.


\subsection{Example: $\pmb{D=-8}$}

This example illustrates Algorithm 1 with the case $d=2$.
The ring $\Z[ \sqrt{-2} ]$ has discriminant $D=-8$.
The elliptic curve
\[
    E : y^2 = x^3 + x^2 - 3x + 1
\]
has $j$-invariant equal to 8000 and its endomorphism ring
is isomorphic to $\Z[ \sqrt{-2} ]$.

We seek a 2-isogeny to a curve with $j$-invariant also equal
to 8000.
Consider the rational
2-isogeny whose kernel is generated by the 2-torsion
point $(1, 0)$.
The equations for this isogeny (found using \cite{velu}) are
\[
   (x,y)  \longmapsto ((3x^2 - 2x + 5) / (3(x - 1)) , y (x^2- 2x - 1) /
    (x - 1)^2 )
\]
and the
image under this isogeny is the elliptic curve
\[
   E' : y^2 = x^3 - 40 x/3 - 448/27.
\]
The curve $E'$ has $j(E') = 8000$ but it is not isomorphic
to $E$ over $\Q$.  There is an isomorphism from $E'$ to $E$ over
$\Q( \sqrt{-2})$ given by
\[
    (x , y ) \longmapsto (-x/2 - 1/3 , \sqrt{-2} y /4 )
\]

The composition of the 2-isogeny and the isomorphism
gives a distortion map $\psi : E \rightarrow E$ which,
by Proposition~\ref{red},
is suitable for our application.
This can be used for $E$ over $\F_p$ whenever $p$
is inert in $\Q( \sqrt{-2})$ (i.e., $p \equiv 5, 7 \mod{8}$).

We note that nicer equations
in this case are known, see
Section 14B of \cite{cox}, \cite{glv} or \cite{rm}.

\subsection{Example: $\pmb{D=-7}$}

We consider the CM curve with $j$-invariant $-3375$
and endomorphism ring $\Z [ (1 + \sqrt{-7})/2 ]$.
The units of this ring  are simply $\pm 1 $.
We consider the curve equation (obtained from Cremona's
tables \cite{cremona})
\[
    E : y^2 + xy = x^3 - x^2 - 2x - 1 .
\]
By Deuring's reduction theorem (see Lang \cite{lang}
Theorem 12 on page 182)
this curve has supersingular reduction modulo $p$
whenever $p=7$ or $(\frac{-7}{p})=-1$ (i.e.,
$p \equiv 2, 5, 6 \mod{7}$).
When $E$ is supersingular modulo $p$ then
$\#E( \F_p ) = p+1$ and the embedding degree is $k=2$.

We seek a non-rational isogeny from $E$ to itself.
Since $\Z [ (1 + \sqrt{-7})/2 ]$ contains
$\sqrt{-7}$ we could apply Algorithm 1 to get a
7-isogeny.  Instead, we note that $\Z [ (1 + \sqrt{-7})/2 ]$
contains elements of
norm 2 and so we should be able to find a 2-isogeny.

Since the kernel of a 2-isogeny is an element of
order 2, we start by finding the 2-torsion on $E$ in
characteristic zero.
Recall that a point $P=(x,y)$ has order 2 if
$P = -P$ and in this case $-P = (x, -y-x)$
hence we require that $x = -2y$.
One easily checks that
\[
    E[2] =  \{ 0_E,  (2, -1), (-2\alpha, \alpha),
    (-2 \overline{\alpha}, \overline{\alpha} ) \}
\]
where $\alpha = (5 + \sqrt{-7})/16$.

The isogeny coming from $(2,-1)$ is rational, and turns
out not to be useful.
Hence we apply V\'elu's formulae \cite{velu} to construct an
isogeny with kernel generated by the point
$(-2\alpha, \alpha)$.
Summarising the results, let
\[
   A_4 = (-29 - 105 \sqrt{-7} )/32 \qquad \mbox{and} \qquad
   A_6 = (-849 + 595 \sqrt{-7} )/128
\]
and  define
\begin{eqnarray*}
  X &=& x + (-7 + 21 \sqrt{-7} )/(32x + 20 + 4 \sqrt{-7} ) \\
  Y &=& y - (-7 + 21 \sqrt{-7} )(2x + 2y + (5 +  \sqrt{-7})/8)
    /(8x + 5+  \sqrt{-7} )^2,
\end{eqnarray*}
Then the map $\psi_1(x,y) = (X, Y)$ is a 2-isogeny from
$E$ to
\[
   E' : Y^2  + XY = X^3 - X^2 + A_4 X + A_6
\]
where $j(E') = -3375$ too.

It remains to compute an isomorphism from $E'$ to $E$.
Let
\begin{eqnarray*}
u &=& (-1 - \sqrt{-7} )/4 \\
r &=& (11 - \sqrt{-7})/32 \\
t &=& (-11 + \sqrt{-7})/64 \\
s &=& (-5 - \sqrt{-7})/8.
\end{eqnarray*}
Then the mapping
$\psi_2 (X, Y) = (u^2 X + r, u^3 Y + u^2sX + t )$
is an isomorphism from $E'$ to $E$.

Defining $\psi(x,y) = \psi_2( \psi_1( x, y ))$ we
obtain our distortion map from $E$ to $E$.
In practice, it is easier to store the isogenies
separately and to compute the distortion map by
computing the composition.

Proposition~\ref{red} does not apply to this
map, so we give a direct proof that it is suitable.
Consider a point $Q$
on the reduction of $E$  over $\F_{p^m}$
($m$ odd)
where $p$ is inert in $\Q( \sqrt{-7} )$.
Let $\pi$ be the $p^m$-power Frobenius.
If $Q \not\in \ker(\psi)$
then we
show that $\pi \psi (Q) \ne \psi \pi(Q)$.
The $x$-coordinate of the composition
of the isogeny and the isomorphism is
\[
   \frac{-3 + \sqrt{-7} }{8} x + \frac{(-63 - 35 \sqrt{-7})/16}{8x + 5
    + \sqrt{-7}} + \frac{11 - \sqrt{-7}}{32}.
\]
Since $\pi$ maps $\sqrt{-7} \in \F_{q^2}$ to $-\sqrt{-7}$ it is clear
that we cannot have $\pi \psi (Q) = \psi \pi(Q)$ for
any point $Q$ except the points in the kernel of $\psi$.


As noted above, this example shows that
Algorithm 1 does not necessarily
provide an endomorphism of minimal degree.
Finally, we note that nicer equations
in this case are known, see
Section 7.2.3 of \cite{cohen}, \cite{glv} or \cite{rm}.

\section{Remaining hard problems}

In the ordinary case, Verheul \cite{verheul} has
shown that there are no distortion maps.
In this case it seems that DDH is hard in both
eigenspaces for the Frobenius map.

To solve the DDH problem in the small field one might try
to invert the trace map.
In fact it is trivial to find pre-images under the
trace map (for example, given $R \in E( \F_q )$ a pre-image
would be $k^{-1}R$)
but it seem to be difficult to find pre-images in a coherent
way without using some kind of non-rational group homomorphism.

It remains an open problem to either show that
DDH is easy on ordinary elliptic curves in all cases, or
to give evidence that the problem is hard in the two
cases remaining (i.e., the two eigenspaces of Frobenius).

The generalisation of these results to the case of
abelian varieties of higher dimension seems to be
hard.  In particular, our algorithm relies on modular
equations to compute isogenies, and it is a well-known
open problem to extend these techniques to the higher-dimensional
case.

\section{Acknowledgements}

We are grateful to
Paulo Barreto, Florian Hess, Takakazu Satoh,
Alice Silverberg and Eric Verheul
for comments on an earlier version of the paper.

\end{document}